\documentclass[12pt]{amsart}
\usepackage{amssymb}
\numberwithin{equation}{section}
\def\T{\text}
\def\1#1{\overline{#1}}
\def\2#1{\widetilde{#1}}
\def\3#1{\widehat{#1}}
\def\4#1{\mathbb{#1}}
\def\5#1{\frak{#1}}
\def\6#1{{\mathcal{#1}}}

\def\C{{\4C}}
\def\R{{\4R}}

\def\Z{{\4Z}}

\begin{document}
\title[CR extension...]{CR extension from manifolds of higher type}
\author[L.~Baracco,  G.~Zampieri]{Luca Baracco, Giuseppe Zampieri}
\address{Dipartimento di Matematica, Universit\`a di Padova, via 
Belzoni 7, 35131 Padova, Italy}
\email{baracco@math.unipd.it,  
zampieri@math.unipd.it}
\maketitle
\def\Label#1{\label{#1}}
\def\cn{{\C^n}}
\def\cnn{{\C^{n'}}}
\def\ocn{\2{\C^n}}
\def\ocnn{\2{\C^{n'}}}
\def\const{{\rm const}}
\def\rk{{\rm rank\,}}
\def\id{{\sf id}}
\def\aut{{\sf aut}}
\def\Aut{{\sf Aut}}
\def\CR{{\rm CR}}
\def\GL{{\sf GL}}
\def\Re{{\sf Re}\,}
\def\Im{{\sf Im}\,}
\def\codim{{\rm codim}}
\def\crd{\dim_{{\rm CR}}}
\def\crc{{\rm codim_{CR}}}

\def\phi{\varphi}
\def\eps{\varepsilon}
\def\d{\partial}
\def\a{\alpha}
\def\b{\beta}
\def\g{\gamma}
\def\G{\Gamma}
\def\D{\Delta}
\def\Om{\Omega}
\def\k{\kappa}
\def\l{\lambda}
\def\L{\Lambda}
\def\z{{\bar z}}
\def\w{{\bar w}}
\def\Z{{\1Z}}
\def\t{{\tau}}
\def\th{\theta}
\emergencystretch15pt
\frenchspacing
\newtheorem{Thm}{Theorem}[section]
\newtheorem{Cor}[Thm]{Corollary}
\newtheorem{Pro}[Thm]{Proposition}
\newtheorem{Lem}[Thm]{Lemma}

\theoremstyle{definition}\newtheorem{Def}[Thm]{Definition}

\theoremstyle{remark}
\newtheorem{Rem}[Thm]{Remark}
\newtheorem{Exa}[Thm]{Example}
\newtheorem{Exs}[Thm]{Examples}

\def\bl{\begin{Lem}}
\def\el{\end{Lem}}
\def\bp{\begin{Pro}}
\def\ep{\end{Pro}}
\def\bt{\begin{Thm}}
\def\et{\end{Thm}}
\def\bc{\begin{Cor}}
\def\ec{\end{Cor}}
\def\bd{\begin{Def}}
\def\ed{\end{Def}}
\def\br{\begin{Rem}}
\def\er{\end{Rem}}
\def\be{\begin{Exa}}
\def\ee{\end{Exa}}
\def\bpf{\begin{proof}}
\def\epf{\end{proof}}
\def\ben{\begin{enumerate}}
\def\een{\end{enumerate}}

\section{notations, generalities, and statements}
\Label{notations}
Let $M$ be a real submanifold of $\C^N$ of codimension $l$ in a neighborhood of a point $p_o$. We assume all through the paper that $M$ is generic which means that its tangent plane $T_{p_o}M$ is not contained in any complex proper subspace of $\C^N$. A wedge $W$ in $\C^N$ is a domain which satisfies, for an open cone $\Gamma$ and a neighborhood $B$ of $p_o$,
\begin{equation}
\Label{1.1}
\left((M\cap B)+\Gamma\right)\cap B\subset V.
\end{equation}
The maximal cone for any whose proper subcone and for suitable $B$   \eqref{1.1} holds is invariant under $T_{p_o}M$ and is therefore identified to a cone $\Gamma$ in the normal space $(T_M\C^N)_{p_o}$, the so called ``directional" cone of $W$ at $p_o$.

We will deal with the space $\T{CR}_M$ of the continuous CR functions on $M$ that is the solutions $f$ of the equation $\bar\partial_Mf=0$ where $\bar\partial_Mf$ denotes the component of $\bar\partial f$ tangential to $M$. (When $f$ is not $C^1$ the equation $\bar\partial_Mf=0$ must be understood in the sense of currents.) A large class of CR functions is described as ``topological" boundary values. Thus, if $F$ is a holomorphic function on a wedge $W$ with edge $M$, continuous up to $M$,  then its boundary value $f=b(F)$ is a CR function on $M$ due to $\bar\partial_Mf=b(\bar\partial F)(=0$). Note that by the {\sc Ajrapetyan-Henkin} edge of the wedge theorem \cite{AH81}, there is a maximal directional cone $\Gamma$ for wedge extendibility of $f=b(F)$. In particular, if we denote by $\Gamma^*$ the polar of this maximal cone, we can meaningfully define the analytic wave front set of $b(F)$ by
$$
WF(b(F))_{p_o}=-\Gamma^*.
$$
The notion of wave front set for CR functions more general than just boundary values requires heavy microlocal machinery \cite{BR87} and goes beyond the purpose of the present presentation. We write now complex coordinates as $(z,w)\in\C^l\times\C^n=\C^N,\,\,z=x+iy$, and suppose that $M$ is defined in a neighborhood of $p_o=0$ by a system of equations $y_j=h_j\,\,j=1,\dots,l$ with $h(0)=0$ and $\partial h(0)=0$; we also write $r=(r_j)_j=(-y_j+h_j)_j$. We select one of the $w$-coordinates, say $w_1$, and define $\tilde M:=M\cap(\C^l_z\times\C^1_{w_1})$. 
We decompose $l$ as $l=l_1+\dots+l_r$, write $I_1=(1,\dots,l_1),\dots,I_r=(\sum_{j\leq r-1}l_j,\dots,l)$ and decompose $z$ as $z=(z_{I_1},\dots,z_{I_r})$. For a set of integers $m_1<\dots<m_r$, where $m_r$ is possibly $+\infty$, we define the notions of ``weighted" homogeneity and vanishing order.  
 For a function $g=g(x_{I_1},\dots,x_{I_j},w_1)$, with $j\leq r$,  we say that $g$ is homogeneous of ``weight" $m_j$ when $h(t^{m_1}x_{I_1},\dots,t^{m_j}x_{I_j},tw_1)$ is a homogeneous polynomial in $t$ of degree $m_j$. We say that $g$ is infinitesimal of weight $m_j$ and write $h=\mathcal O^{m_j}$ when
$$
g(t^{m_1}x_{I_1},\dots,t^{m_j}x_{I_j},tw_1)=O(t^{m_j}).
$$
A special definition is needed for $j=r$ and $m_r=+\infty$. In this case  we say that $g$ is infinitesimal of weight $+\infty$ and write $g=\mathcal O^{+\infty}$ when $g(t^{m_1}x_{I_1},\dots,t^{m}x_{I_r},tw_1)=O(t^m)$ for any $m$. In other terms, $g$ is divisible by some monomial in $x_{I_r}$. 
We recall what equations for $\tilde M$ in ``Bloom-Graham normal form" mean. Intrisically associated to $\tilde M$ there are integers $m_1<\dots<m_r$, the so called ``H\"ormander numbers", and $l_1,\dots,l_r$ with $\sum_jl_j=l$, their respective ``multiplicities". For $m_r<+\infty$, in suitable coordinates at $p_o$, $\tilde M$ is described by equations
\begin{equation}
\Label{1.2}
\begin{cases}    
y_{I_1}=P_{I_1}(w_1)+\mathcal O^{m_1+1},
\\
y_{I_2}=P_{I_2}(x_{I_1},w_1)+\mathcal O^{m_2+1},
\\
\dots,
\\
y_{I_r}=P_{I_r}(x_{I_1},\dots,x_{I_{r-1}},w_1)+\mathcal O^{m_r+1},
\end{cases}
\end{equation}
with each $P_{I_{j}}$ homogeneous of degree $m_j$
and such that for any $\xi^o\in\R^{l_j}$, $\langle\xi^o,P_{I_j}\rangle$ is not $\tilde M$-pluriharmonic.
(A homogeneous polynomial $g$ of weight $m_j$ is said $\tilde M$-pluriharmonic of weight $m_j$ if ther exists $F$ holomorphic in $\C^{l+1}$ such that $g=\Im F|_{\tilde M}+\mathcal O^{m_j+1}$.)
When $m_r=+\infty$, then for any $m$, there are coordinates such that \eqref{1.2} holds with the last equation replaced by $y_{I_r}=\mathcal O^m$.
$\tilde M$ is said to be of ``finite type" when $m_r<+\infty$. $\tilde M$ is said ``semirigid" when each $P_{I_j}$ is a function of $w_1$ only. The similar notions of finite type and semirigidity for $M$ instead of $\tilde M$ apply when one deals with equations of type \eqref{1.2} involving all $w$-variables instead of $w_1$ only. We will see in \S 3 that finite type can be characterized by means of brackets instead of normal equations: iterated commutators of vector fields tangential to $M$, of $(1,0)$ and $(0,1)$ type, up to a certain finite number,  the highest H\"ormander number $m_r$, span the whole complexified tangent bundle $\C\otimes_\R TM$. Let us recall that when $M$ is of finite type, then according to {\sc Tumanov} \cite{T90}, CR functions $f$ are boundary values $f=b(F)$ of holomorphic functions $F$ on a wedge $W$; in particular, in this situation, the notion of wave front set applies to any $f$. 
\bt
\Label{t1.1}
Let $M$ be a generic manifold of $\C^N$ of finite type, and, for a choice of a complex tangent direction $w_1$, let \eqref{1.2} be a normal system of equations for $\tilde M=M\cap(\C^l_z\times\C^1_{w_1}$). We assume that for some $j$, for $\xi^o\in\R^{l_j+\dots+l_r}$ and with the notation $P:=\langle\xi^o,P_{I_j}\rangle$ we have
\begin{equation}
\Label{1.3}
\begin{cases}
P=P(w_1)\T{ for a homogeneous polynomial 
$P$ of degree $m_j$},
\\
P(w_1)\geq0\T{ for $w_1$ in a sector $\mathcal S$  of width $>\frac 
\pi{m_{j}}$}.
\end{cases}
\end{equation}
Then
\begin{equation}
\Label{1.4}
\xi^o\notin \T{WF}(f)\quad\forall f\in CR_M.
\end{equation}.
\et
The proof will follow in \S 2. The first of \eqref{1.3} is  a sort of semirigidity in direction $w_1$ and codirection $\xi^o$. We will exhibit in \S 4 (Proposition~\ref{p8} and Corollary~\ref{c4.1}) a large class of hypersurfaces $M$ for which when \eqref{1.3} is violated, we can find a ``barrier" that is a holomorphic function $F$ with $M\subset\{\Im F<0\}$. In particular, for these $M$, there always exist CR functions $f\in CR_M$ such that $\xi^o\in WF(f)$ for $\xi^o=d(\Im F)$. This shows that the statement in Theorem~\ref{1.1} is sharp.
\br
When $j=1$ the first of \eqref{1.3} is automatically fulfilled. Also, since we are assuming that $P_{I_1}$ is not $\tilde M$-harmonic, then it is divisible by $|w_1|^2$ and therefore it has at most $2m_1-2$ zeroes on the unit circle $|w_1|=1$. In particular, for either of $\pm P$, the second of \eqref{1.3} is satisfied. 
\er
\br
There is a sort of ``hierarchy" between the H\"ormander numbers $m_j$ whose geometric meaning will be fully clear from the proof in \S 2. According to it, \eqref{1.3} for $j>1$, gives the control not of the whole $\T{WF}(f)$ but only of its section $\T{WF}(f)\cap\left(\{0\}\times\dots\times\{0\}\times\R^{l_j+\dots+l_r}\right)$. In fact, the proof of the theorem will consist in proving CR extension in some extra direction $v$ close to the component normal to $M$ of the disc attached to $M$ over $\mathcal S$, and \eqref{1.3} does not give informations for $v$ itself but for  $v_{I_j,\dots,I_r}$.
\er
When $M$ is of finite type and semirigid (in the complex of its arguments $w$), then our proof provides an alternative proof of the extension of any $f$ to a wedge $W$. The first conclusion in this direction is due to \cite{BR87} where a description of $W$ is also given. We improve this description by specifying the vanishing order in a precise direction $w_1$. Also,  the semirigidity in the first of \eqref{1.3} can be released,
as well the hypothesis that the equations are in canonical form as in \eqref{1.2}. What is indeed essential is the weighted vanishing order; non-$\tilde M$-harmonicity in the homogeneous terms is not needed.
 Thus, suppose that $M$ is of finite type and that $\tilde M=M\cap(\C^l\times\C)$ has equations in the (not necessarily normal) form $y_{I_j}=h_{I_j}$ with $h_{I_j}=\mathcal O^{m_j},\,\,j=1,...,r$. 
\bt
\Label{t1.2}
In the above situation suppose that for  $j\leq r$ with $m_j<+\infty$ and for some $\xi^o\in\R^{l_j+\dots+l_r}$, we have for a suitable constant $c$
\begin{equation}
\Label{1.5}
\langle\xi^o,h_{I_j}\rangle>0\T{ for $w_1$ in a sector $\mathcal S$ of width $>\frac{\pi}{m_j}$ and for $|x_{I_i}|<c|w_1|^{m_i}$}.
\end{equation}
Then 
$$
\xi^o\notin WF(f)\quad\forall f\in CR_M.
$$
\et
(If $m_r=+\infty$ in \eqref{1.5} the condition $|x_{I_r}|<|w_1|^{m_r}$ means $|x_{I_r}|<|w_1|^m\,\,\forall m$.) The proof follows in \S 2. 
If $h_{I_j}=P_{I_j}+\mathcal O^{m_j+1}$, then clearly 
the components of $h_{I_j}$ have the same sign as those of $P_{I_j}$ under the constrain $|x_{I_i}|\leq c|w_1|^{m_i}$; hence the second of \eqref{1.3} implies \eqref{1.5}. This shows that Theorem~\ref{t1.1} is a particular case of Theorem~\ref{t1.2}.

There are two main streams of CR extension: umprecised extension 
through {\em minimality}; extension in Levi or higher type 
directions. As for the first, it was completely solved by {\sc Tumanov} in \cite{T90} (cf. also {\sc Trepreau} \cite{Tr85} in case $M$ is a hypersurface). He introduced the notion of ``minimality" of $M$ as the absence of proper submanifolds $S\subset M$ with the same CR structure as $M$ that is $T^\C S=T^\C M|_S$. Note that when $M$ is of finite type then it is minimal. (First, finite type in the sense that $m_r<+\infty$ in a system of normal equations, is equivalent to ``finite bracket type" according to the subsequent discussion of \S 3. But then the presence of $S$ as above would force all brackets to belong to $\C\otimes_\R TS$.) He then proved that if $M$ is minimal, then there exist arbitrarily small discs of ``defect" $0$ and hence endowed with infinitesimal deformations which span all normal directions to $M$. Collecting all these directions by the edge of the wedge theorem of \cite{AH81} he got a common wedge $W$ to which all CR functions are forced to extend. Necessity of minimality for such an extension is on its hand a simpler result (and even trivial if $M$ is a hypersurface). However, a precise description of $W$ has not yet been found. Our paper aims at this attempt and deals with extension in directions produced by higher type commutators.  Let us briefly recall the related literature. The first theorems go back 
to {\sc Ajrapetyan-Henkin} \cite{AH81} and {\sc Boggess-Polking} 
\cite{BP82} and state extension in directions of the Levi cone. Next 
{\sc Boggess-Pitts} \cite{BP85} proved extension in directions shown 
by iterated brackets up to the first H\"ormander number. More 
recently the authors obtained in collaboration with {\sc Zaitsev}
generalizations to the case of CR functions not defined on the whole 
$M$ but, instead, on a subwedge $V\subset M$. Let us point out the 
main novelties of the present paper. In the equations  \eqref{1.2} 
the weighted homogeneity degrees $m_1<m_2<...$ are calculated with 
respect to $w_1$ and not to the complex of the variables. Also, $m_{j}$ 
is not the smallest among the $m_i$'s. 
On the contrary, most of other classical CR extension criteria 
concern the first H\"ormander number: in all $w$ directions as in 
\cite{BP85}, or at least in one selected direction as in \cite{BR87}  Th. 11. 
(Let us point out that it seems that the method of \cite{BP85} can be 
adapted to treat also this second situation though this is absent 
from their statements.)
This  paper \cite{BR87}, whose conclusions are the closest to ours, 
gives indeed extension also related to further H\"ormander numbers, 
as e.g. in Theorem 8. But in this case its method founded on Fourier 
calculus, requires an assumption of semirigidity in the complex of 
the equations and of the arguments $w$. To explain the difference, 
let's consider for instance 
in $\C^4$ with coordinates $(z_1,z_2,w_1,w_2)$ the 
manifold $M$ defined by 
\begin{equation}\Label{mainexample}
\begin{cases}
y_1=|w_1|^2+|w_2|^2+f_1(x_1,x_2,w_1,w_2)
\\
y_2=|w_2|^4+x_1|w_1|^2+f_2(x_1,x_2,w_1,w_2)
\end{cases}
\end{equation}
where 
$f_j=\mathcal O^{2j+1}$.
Extension in direction $v^1\sim(1,0)$ is clear according to all 
authors. For extension in directions with non-trivial 
$y_2$-component, we notice that the method by \cite{BR87} fails because 
of the lack of semirigidity. 
(Also, \cite{BP85} and \cite{BR87} Th.11 cannot be applied because 
$w_2$ appears in the second equation in a higher homogeneity degree 
than in the first.)
However our Theorem~\ref{t1.1} applies for sectors in the $w_2$-plane, 
and yields extension in direction 
$v^2=(1+\sigma(\eta),\eta^2+o(\eta^2))$ with $\sigma(\eta)$ infinitesimal with $\eta$. 
Our generalization goes also in another direction (though this was 
already achieved in \cite{BZZ02}). We 
are able to 
obtain extension in more general situations and to a larger set of 
directions. Thus, for instance, let $M$ be the manifold in $\C^3$ 
defined by $(y_1=|w|^4+a|w|^2\Re w^2\,,\,y_2=|w|^4)$. \cite{BP85} 
gives extension for $a>2$ in directions 
which satisfy $y_1>-|y_2|(\frac a2-1)$. On  the other hand, by the proof of our 
Theorems~\ref{t1.1} and \ref{t1.2}, we have extension when $a>\sqrt2$ in directions 
satisfying $y_1>-|y_2|(\frac a{\sqrt2}-1)$. We refer to the 
subsequent \S 4 for a complete proof of these 
claims.
\newline
{\sc Aknowledgments.} The authors are grateful to Professor Dmitri 
Zaitsev for many valuable discussions. In particular they owe to him 
the idea of putting the proof of Theorem~\ref{Fgamma} in the elegant 
form of Lemma~\ref{Cgamma}.

\section{Proof of Theorems~\ref{t1.1}, \ref{t1.2}.} 

 {\bf (a) Preliminaries on $\mathcal F^\alpha$ spaces.} Let 
$0<\alpha<1$ and denote by $\tau=re^{\T{i}\theta}$ the variable in 
the standard disc $\Delta$.
Let us recall from 
\cite{T96}, \cite{ZZ01} and \cite{ZZ02} some basics about 
attaching analytic 
discs to $M$ in the subclasses $\mathcal F^\alpha$ of the H\"older 
classes $C^\alpha$.
These are the spaces of real continuous functions 
$\sigma(\theta),\,\,\theta\in[-\pi,\pi]$ which are $C^{1,\alpha}$ out 
of $0$ and for which the following norm is finite
$$
||\sigma||_{\mathcal 
F^\alpha}:=||\sigma||_{C^0}+||\theta\sigma^{(1)}||_{C^\alpha}.
$$
(Here $\cdot^{(1)}$ denotes the first derivative.) We remark that for 
$\sigma\in\mathcal F^\alpha$ we 
must have $\theta\sigma^{(1)}|_{\theta=0}=0$ for otherwise 
$\theta\sigma^{(1)}\to c\neq0$ which implies $|\sigma|\geq 
\T{log}\,\frac{|c|}2+\T{log}\,|\theta|$ which contradicts the 
boundedness of $\sigma$.
This shows that $\mathcal F^\alpha$ is continuously embedded into $C^\alpha$. It is easy to check that $\mathcal F^\alpha$ is a Banach algebra. Also, if $\sigma_i\in \mathcal F^{\alpha_i}\,\,i=1,2$, then $\sigma_1\cdot\sigma_2\in\mathcal F^{\alpha_1+\alpha_2}$ for $\mu_1+\mu_2<1$, resp. $\sigma_1\cdot\sigma_2\in C^{1,\beta}$ with $\beta:=(\mu_1+\mu_2)-1$ for $\mu_1+\mu_2>1$. In both cases the multiplication is continuous with values in the respective spaces.

Let $T_1$ denote the Hilbert transform normalized by the condition 
$T_1(\cdot)(1)=0$;
it is easy to see that  $T_1$ is a bounded operator in $\mathcal 
F^\alpha$. We come back to our manifold $M$. We write coordinates in 
$\C^N\simeq \C^l\times\C^{N-l}$ as $(z,w)$ with $z=x+iy$, 
choose a distinguished direction, say $w_1$, 
and 
describe 
$\tilde M:=M\cap(\C^l_z\times\C^1_{w_1})$ by the system of equations $y_{I_i}=h_{I_i}(x,w)$ with $h_{I_i}=\mathcal O^{m_i}$. 
(The Bloom-Graham normal form is not needed.) 
We will consider in $\C^N$ analytic discs 
$A(\tau)=(z(\tau),w(\tau))\,\,\tau\in\Delta$ (the standard disc in 
$\C$), {\em attached} to $\tilde M$ that is satisfying 
$A(\partial\Delta)\subset \tilde M$. If we prescribe an analytic function 
$w_1(\tau)\,\,\tau\in\Delta$, the so called {\em CR component}, and a 
point $p=(z,w_1)$ with $y=h(x,w_1)$, and look for an analytic 
completion $z(\tau)$ for $A(\tau)=(z(\tau),w_1(\tau))$ with $A(1)=p$, 
we are lead to the {\em Bishop's equation}
\begin{equation}\Label{Bishop}
u(\tau)=-T_1h(u(\tau)+x,w_1(\tau))\quad \tau\in\partial\Delta.
\end{equation}
In fact if $u(\tau)$ solves \ref{Bishop}, then if set 
$z(\tau)=u(\tau)+iv(\tau)+z$, we obtain that 
$A(\tau)=(z(\tau),w_1(\tau))$ is 
holomorphic, 
$v(\tau)=h(u(\tau),w_1(\tau))\,\,\forall\tau\in\partial\Delta$ and 
finally $A(1)=p$. We will consider equation \eqref{Bishop} in the 
spaces 
$\mathcal F^{\alpha}$, $\mathcal F^{m_i\alpha}$  and $C^{1,\beta}$ for which $T_1$ is bounded. 
We will also use the composition properties of $h_{I_i}$ for $i\geq j$ with functions in 
the above classes as stated in \cite{BZZ02}. To take advantage of this 
composition we will assume $m_j\alpha>1$ (and, to be sharp, 
$\alpha(m_{j}-1)<1$). Here is our main technical tool.

\bp\Label{p4}
Let $h_{I_i}$ be of class $C^{m_i+3}$,
and satisfy $h_{I_i}=\mathcal O^{m_i}$. Then for any $\epsilon$ there 
is $\delta$ 
such that if $||h_{I_i}||_{C^{1,\alpha}}<\delta,\,\,||w_1||_{\mathcal 
F^\alpha}<\delta,\,\,|x|<\delta$, then the equation \eqref{Bishop} 
has an unique 
solution $u\in\mathcal F^\alpha$ with $||u\|_{\mathcal 
F^\alpha}<\epsilon$. 
Moreover, $u_{I_1}\in\mathcal F^{m_1\alpha},\dots,u_{I_{j-1}}\in\mathcal F^{m_{j-1}\alpha}$ and $(u_{I_j},\dots,u_{I_r})\in C^{1,\beta}$ for $\beta=m_j\alpha -1$ and, if $w_1$ 
depends on some 
parameter $\lambda\in\R^d$ so that $\lambda\mapsto 
{w_1}_\lambda,\,\,\R^d\to\mathcal F^\alpha$ is $C^k$ for $k\leq m_i$, 
then also $\lambda,x\mapsto (u_{I_i})_{\lambda \,x},\,\,\R^{d+l}\to 
C^{1,\beta}$ 
is $C^k$. 
In particular there exist mixed derivatives in $\lambda, x$ and $r$ 
up to order $k$ and $1$ respectively, and they commute that is

\begin{equation}\Label{7}
 \partial_r\partial^{k'}_{\lambda\,x}u= \partial_{\lambda\, 
x}^{k'}\partial_r 
u\,\,\forall k'\leq k.
\end{equation}
\ep

\bpf
One first solves Bishop's equation \eqref{Bishop} in the $\mathcal 
F^\alpha$ - 
spaces by the aid of the implicit function theorem. To this end one 
considers 
the mapping $F\colon (\lambda,x,w_1,u)\mapsto 
u-T_1h(u+x,w_1),\,\,\R^d\times\R^l\times\mathcal F^\alpha\times\mathcal 
F^\alpha\to \mathcal 
F^\alpha$. Then for the partial Jacobian $\partial_uF$ with respect 
to $u$, one has $\partial_uF:\,\dot u\mapsto\dot u-T_1\partial_xh\dot 
u$. In 
particular if we evaluate at $(\lambda,x,w_1,u)=(0,0,0,0)$, then this 
is invertible since $\partial_x h|_0=0$. The differentiability with 
respect to the parameters in the space $\mathcal F^\alpha$ is also clear 
in view of \cite[Prop. 11]{ZZ01}.

We show that the components $u_{I_i}\,\,i\geq j$ of the solution to Bishop's equation, as well as their 
harmonic conjugates $v_{I_i}$, are in fact in $\mathcal F^{m_i\alpha}$ for $i<j$ (resp. $C^{1,\beta}$ for $i\geq j$)  with $\beta:=m_j\alpha-1$
and also prove differentiability in the parameters with values in this 
space.
 The key point is that the composition $\phi((1-\tau)^\alpha)$, and 
in bigger 
 generality $\phi(w_1)$ for $w_1\in \mathcal F^\alpha,\,\,w_1(1)=0$, with 
 $\phi=O^{m_i}$ belongs to $\mathcal F^{m_i\alpha}$ for $i<j$ (resp. $C^{1,\beta}$ for $i\geq j$). We 
 put $z(\tau)=u(\tau)+iv(\tau)+z$ with $v=T_1 u$ and $z=x+ih(x,w)$ 
and also write $\tau=e^{i\theta}$ on $\partial\Delta$. We can check 
that if $z_{I_i}(\tau)\in\mathcal F^{k\alpha}$ for $k\leq m_i-2$, then in 
fact $z(\tau)\in \mathcal F^{(k+1)\alpha}$. In fact $v$ gains regularity at each step because  $h_{I_i}=\mathcal O^{m_i}$
together with the fact that if $\sigma(\theta)\in 
\mathcal 
F^{k\alpha}$ and $\sigma(0)=0$, then 
$|\theta|^\alpha\sigma(\theta)\in\mathcal F^{(k+1)\alpha}$ due to 
\begin{equation}
\begin{split}
|(|\theta|^\alpha\sigma(\theta))^{(1)}|&=||\theta|^{\alpha-1}\sigma(\theta)+
|\theta|^\alpha\sigma^{(1)}(\theta)|
\\
&\leq c|\theta|^{(k+1)\alpha-1}.
\end{split}
\end{equation}
But the Hilbert transform interchanges the 
$\mathcal F^{(k+1)\alpha}$ regularity from $v$ to $u$ and 
thus $z(e^{i\theta})\in\mathcal F^{(k+1)\alpha}$. This completes the proof when $i<j$. On the other hand, when $i\geq j$, 
 in order to pass from $\mathcal F^{(m_i-1)\alpha}$ 
to $C^{1,\beta}$, 
we have to prove that $(\theta^\alpha u)^{(1)}=\theta^{\alpha-1}u+\theta^{\alpha-1}(\theta u^{(1)})$ belongs to $C^\beta$. But in fact, since both $u$ and $\theta u^{(1)}$ are in $C^{(m_i-1)\alpha}$ and are $0$ at $\theta=0$, then their product by $\theta^{\alpha-1}$ is in $C^\beta$ as one can easily check by the Hardy-Littlewood principle.
 It follows that $(\theta^\alpha u)^{(1)}\in C^\beta$ and hence $\theta^\alpha u\in C^{1,\beta}$. 
Thus  $u(e^{i\theta})$ and hence $z(e^{i\theta})$ itself is in $C^{1,\beta}$. 
As for the differentiability on $x$ and on the parameters, 
it is a variant  of [Proposition 15] in \cite{BZZ02} 
by the same {\em feed back} argument as above.

\epf

We can think of the family of discs produced by the above statement 
as a {\em deformation} of the disc $A(\tau)\equiv0$ which is a 
trivial solution to Bishop's equation. By the next statement we will 
show how it is possible to make infinitesimal deformations of discs 
which are no more small. 
\bp
\Label{p5}
Let $h_{I_i}\in C^{m_i+3}$ satisfy $h_{I_i}=\mathcal O^{m_i}$, let $\tilde w_1(\tau)\in 
C^{1,\beta},\,\,\tilde w_1(1)=0$ be small in $\mathcal F^\alpha$ (not 
necessarily in $C^{1,\beta}$), and let $\tilde u(\tau)\in \mathcal F^\alpha$ be 
a solution of Bishop's equation $\tilde u=-T_1h(\tilde u,\tilde w)$; in particular $\tilde u_{I_i}\in C^{1,\beta}$ for any $i\geq j$ according to Proposition~\ref{p4}.  
Then for any 
$w_1(\tau)$ with $||w_1-\tilde w_1||_{C^{1,\beta}}<\delta,\,\,|x|<\delta$ 
there 
is an unique solution $u\in\mathcal F^\alpha$ with $u_{I_i}(\tau)\in C^{1,\beta}\,\,\forall i\geq j$ of Bishop's equation 
with $||u_{I_i}-\tilde u_{I_i}||_{C^{1,\beta}}<\epsilon\,\,\forall i\leq j$. Moreover, if 
$\lambda\mapsto 
(w_1)_\lambda$ is $C^k,\,\,k\leq m_i$, then also $\lambda,x\mapsto 
(u_{I_i})_{\lambda}$ is $C^k$.
\ep
\bpf
In the present situation we define $F\colon \R^d\times\R^l\times 
C^{1,\beta}$ similarly as in the proof of Proposition \ref{p4} and 
wish to 
prove that $\partial_u F$ is still invertible. For this purpose it is 
enough to show that 
$\partial_uh_{I_i}(\tilde u,\tilde w)$  is small in $C^{1,\beta}$ - norm. 
But in fact 
recall that $|\partial_xh_{I_i}(u,w)|=O(|w|^2)$ and therefore
$||\partial_xh_{I_i}(\tilde u,\tilde w)^{(1)}||_{C^\beta}\leq 
c||\tilde w_1||_{C^\beta}||\tilde w_1^{(1)}||_{C^\beta}\leq\epsilon$.

\epf

{\bf (b) Construction of a singular disc attached to $M$ 
with controlled normal component.} Let us suppose that \eqref{1.5} be fulfilled. It is not restrictive to assume that the 
sector $\C_{w_1}$ where 
$g\geq 0$ contains $(1-\tau)^\alpha ie_{l+1},\,\,\tau\in\Delta$. 
(Here  $e_{l+1}$ is the unit vector of the $w_1$
- plane.) 
Let  $\alpha$ satisfy $\alpha 
m_j>1,\,\,\alpha(m_j-1)<1$.
We  define, for a small real parameter $\eta>0$:
\begin{equation}\label{9}
w_1(\tau)=(w_{1})_\eta(\tau):=\eta(1-\tau)^\alpha ie_{l+1}.
\end{equation}
We attach to $\tilde M$ a family of $\mathcal F^\alpha$ - discs 
$A(\tau)=A_\eta(\tau)$ whose "$w_1$-component" is $w_1(\tau)$. 
We recall from (a) that for $i\geq j$ we have $\eta\mapsto (z_{I_i})_\eta(\tau),\,\,\R\to 
C^{1,\beta}$ 
is $C^{m_i}$. We also write $z_{I_i}(\tau)$ instead of $(z_{I_i})_\eta(\tau)$, 
$z_{I_i}(\tau)=u_{I_i}(\tau)+iT_1v_{I_i}(\tau)$, and finally $A(\tau)=(z(\tau),w(\tau))$.
 We note that we have
\begin{equation}\Label{10}
\partial_\eta^sv_{I_i}|_{\eta=0}\equiv 0,\quad \partial ^s_\eta 
u_{I_i}|_{\eta=0}\equiv0\quad \forall s\leq m_i-1.
\end{equation}
This is clear for $s=0,1$. If it is true for any 
$s\leq m_i-2$, then it is also true for $s=m_i-1$ due to $h_{I_i}=\mathcal O^{m_i}$ by a ``feed-back" procedure. If we then Taylor-expand $\partial_rv_{I_i}$ at $\eta=0$, we get
\begin{equation}\Label{11}
\partial_rv_{I_i}=\frac{\partial^{m_i}_\eta\partial_rv_{I_i}|_{\eta=0}}{m_i!}\eta^{m_i}+o(\eta^{m_i}).
\end{equation}
By a similar argument we can also prove that
\begin{equation}
\Label{feedback}
|v_{I_i}|\leq c|w_1|^{m_i},\quad |u_{I_i}|\leq c|w_1|^{m_i}.
\end{equation}
In fact, in the classes $\mathcal F^{k\alpha}$ regularity and vanishing order are coincident: thus the equation $v_{I_i}=h_{I_i}$ gives control of the vanishing order of $v_{I_i}$ which is transferred as regularity to $u_{I_i}$ through Hilbert transform, and again as vanishing order to $v_{I_i}$. In this way we can prove that each $v_{I_i}$ and $u_{I_i}$ belongs to $\mathcal F^{m_i\alpha}$ (and also to $C^{[m_i\alpha],\{m_i\alpha\}}$ where $[m_i\alpha]$, resp. $\{m_i\alpha\}$, is the integer, resp. fractional, part of $m_i\alpha$. 
Recall that if $\xi_o$ is, say, the unit vector in the $l':=l_{1+...+l_{j-1}+1}$-direction, we 
have $\langle\xi_o,h\rangle\geq0$ if $w_1$ is in a sector $\mathcal S$ of width 
$>\frac 
{m_j}{\pi})$ and if $|x_{I_i}|\leq c|w_1|^{m_i}$. We first observe that this latter condition $|x_{I_i}|\leq c|w_1|^{m_i}$ is automatically fulfilled by the components $x_{I_i}=u_{I_i}$ of our discs $A(\tau)$ due to \eqref{feedback}. We  show now that  $\partial_rv_{l'}<0$. In fact
we have in this situation
\begin{equation}\Label{13}
\langle\xi_o,\partial^{m_j}_\eta v_{l'}\rangle|_{\eta=0}\geq0 \forall \tau\in\bar\Delta.
\end{equation}
 Hence \eqref{1.5} yields, 
through Hopf's Lemma 
\begin{equation}\Label{14}
\langle\xi_o,\partial_r\partial^{m_j}_\eta v_{l'}\rangle|_{\tau=1\,\eta=0}=-c<0.
\end{equation}
By \eqref{11} we conclude $\langle\xi_o,
\partial_rv_{l'}\rangle|_{\tau=1}=-c'\eta^{m_j}<0,$ 
for any $\eta$ sufficiently small. We fix such a small $\eta$ and, by 
rescaling, we even suppose $\eta=1$ and
define $v_o=\partial_rv|_{\tau=1}$. According to \eqref{14} we have 
$\langle v_o,\xi_o\rangle<-\frac c2$.

{\bf (c) Polynomial approximation of $(1-\tau)^\alpha$ in $\mathcal 
F^\gamma(\bar\Delta)$ for $\gamma<\alpha$.}
We have the Taylor expansion
\begin{equation}
\label{c1}
\begin{split}
(1-\tau)^\alpha&=1-\alpha\tau-\frac{\alpha(1-\alpha)}{2!}
\tau^2-\frac{\alpha(1-\alpha)(2-\alpha)}{3!}\tau^3+\dots
\\
&=1-\underset{n=1}{\overset{+\infty}\sum}\left|{\alpha \choose 
n}\right|\tau^n.
\end{split}
\end{equation}
We call 
$S_N=S_N(\tau)$ its the partial sum of the series \eqref{c1} for 
$1\leq n\leq N$.
Our goal is  to prove next
\bt
\Label{Fgamma}
We have
\begin{equation}
\label{c2}
S_N(\tau)\to(1-\tau)^\alpha\text{ in $\mathcal F^\gamma(\bar\Delta)$ 
for any $\gamma<\alpha$}.
\end{equation}
\et
Before giving the proof of Theorem~\ref{Fgamma}, let us recall that 
$\|\sigma\|_{\mathcal 
F^\gamma}=\|\sigma\|_{C^0}+\|(1-\tau)\sigma'\|_{C^\gamma}$.
Hence we have to prove that
\begin{gather}
\label{c3}
S_N\to (1-\tau)^\alpha\T{ in } C^0(\bar\Delta),\quad 
(1-\tau){S'_N}
\to-\alpha(1-\tau)^\alpha\text{ in $C^\gamma(\bar\Delta)$}.
\end{gather}
To prove the first of \eqref{c3}  we note that 
 since 
\begin{equation}
\label{c5}
|S'_N(\tau)|\leq\sum_{n=1}^N\left|{\alpha\choose n}
\right|n|\tau|^{n-1}\to\alpha(1-|\tau|)^{\alpha-1},
\end{equation} 
then in particular the partial sums $|S'_N(\tau)|$ 
are bounded on $\Delta$, uniformly over $N$,
 by $\alpha(1-|\tau|)^{\alpha-1}$. In particular the sequence of 
the $S_N$'s is uniformly continuous in $\bar\Delta$, 
which yields at once the first of \eqref{c3}.
As for the second of \eqref{c3} we note that
\begin{equation}
\label{c6}
|S''_N|\leq\sum_n\left|{\alpha\choose 
n}\right|n(n-1)|\tau|^{n-2}\to\alpha|\alpha-1|(1-|\tau|)^{\alpha-2}.
\end{equation}
It follows
\begin{equation}
\label{c7}
\begin{split}
\left|\left((1-\tau)S'_N\right)'\right|&\leq|S'_N|+|1-\tau||S''_N|
\\
&\leq\alpha(1-|\tau|)^{\alpha-1}+\alpha|\alpha-1|(1-|\tau|)^{\alpha-1}=c(1-|\tau|)^{\alpha-1}.
\end{split}
\end{equation}
To conclude the proof of Theorem~\ref{Fgamma} it suffices to use the 
following one real variable Lemma
\bl
\label{Cgamma}
Let $\{f_N\}$ be a sequence of real functions such that
\begin{equation}
\label{c8}
f_N\to0\T{ in $C^0([0,1-\epsilon])$ for any $\epsilon$},
\end{equation}
and
\begin{equation}
\label{c9}
|f'_N|\leq c(1-t)^{\alpha-1}\T{ in $[0,1)$},
\end{equation}
with $c$ independent of $\epsilon$. Then
\begin{equation}
\label{c10}
f_N\to0 \T{ in $C^\gamma([0,1])$ for any $\gamma<\alpha$}.
\end{equation}
\el

\bpf
We have by integration
$$
|f_N(x)-f_N(y)|\leq c|x-y|^\alpha \T{ (for a different $c$)}.
$$
It follows that for any $\epsilon$ and for suitable 
$\delta=\delta_\epsilon$ we have, when $|x-y|<\delta$
$$
\frac{|f_N(x)-f_N(y)|}{|x-y|^\alpha}|x-y|^{\alpha-\gamma}\leq 
|x-y|^{\alpha-\gamma}<\epsilon.
$$
On the other hand, when $|x-y|\geq\delta$, then
\begin{equation}
\frac{|f_N(x)-f_N(y)|}{|x-y|^\gamma} 
\leq \delta^{-\gamma}|f_N(x)-f_N(y)|\leq
 \delta^{-\gamma}(|f_N(x)|+|f_N(y)|).
\end{equation}
Hence it suffices to prove that $f_N\to0$ in $C^0([0,1])$. 
By \eqref{c9} $\{f_N\}$ is  equicontinuous. 
Given $\epsilon$ we thus have $|f_N(x)-f_N(y)|\leq\epsilon$, 
uniformly on $N$, for any $\xi$ such that $|x-y|\leq\delta$, 
in addition to $\underset{[0,1-\delta]}{\T{sup}}|f_N|<\epsilon$ for 
any $N\geq N_\epsilon$. In conclusion given $x$ we take
 $\xi\in[0,1-\delta]$ with $|x-\xi|<\delta$, and then get, 
for any $N\geq N_\epsilon$
$$
|f_N(x)|\leq|f_N(x)-f_N(\xi)|+|f_N(\xi)|<\epsilon.
$$
This concludes the proof of the Lemma. The proof of 
Theorem~\ref{Fgamma} is also complete.

\epf

{\bf (d) Construction of a smooth disc transversal to $M$ and of its 
infinitesimal deformation.}

We put $w_N(\tau)=S^\alpha_N(\tau)-S_N(1)$, let $u_N$ be the solution 
in $\mathcal F^\gamma$ to Bishop's equation $u_N=-T_1h(u_N,w_N)$, and 
let $z_N=u_N+iv_N$ for $v_N=T_1u_N$. Let $u$ be the solution to 
$u=-T_1h(u,(1-\tau)^\alpha)$, and set $z=u+iv$ for $v=T_1u$. Since
$$
w_N(\tau)\to i(1-\tau)^\alpha\text{ in $\mathcal 
F^\gamma(\bar\Delta)$},
$$
and since $w_N(1)\equiv0\,\forall N$, then $(z_{I_j,...,I_r})_N(\tau)\to z_{I_j,...,I_r}(\tau)$ in 
$C^{1,\beta'}(\bar\Delta)$ by Proposition~\ref{p4}. (Clearly we are 
supposing $\gamma$ close enough to $\alpha$ so that 
$\beta':=m_j\gamma-1>0$.) In particular for any $\epsilon$ and for 
large $N$ the discs $A_N=(z_N,w_N)$ are in $C^{1,\beta'}$ and satisfy
\begin{equation}
\label{d1}
\partial_r(v_{I_j,...,I_r})_N(1))=v'_o\text{ for $|v'_o-v_o|<\epsilon$},
\end{equation}
uniformly in $N$. We call $\tilde A=(\tilde z,\tilde w)$ one of these 
discs.
We are ready to construct a half-space $M_1^+$ in a manifold 
$M_1$ which contains $M$ and gains one more direction by a deformation of the disc $\tilde A$ to which CR functions extend.
For this we consider the Bishop's equation
\begin{equation}
\label{d2}
u=-T_1h(u+x,w+\tilde w),
\end{equation}
for $x\in\R^l,\,w\in\C^{n}$ with $|x|<\delta,\,|w|<\delta$. According 
to Proposition~\ref{p5}, for any $\epsilon$ and for suitable 
$\delta=\delta_\epsilon$ there is an unique solution $u$ which 
satisfies $\|u-\tilde u\|_{C^{1,\beta'}}<\epsilon$ for 
$\beta'<\beta:=k\alpha-1$. We write $p=x+ih(x,w),w)$ with $v=T_1u$, 
and define $A_p(\tau)=p+(u(\tau)+iv(\tau),\tilde w(\tau)$. We also 
write $I_p=A_p|_{[-1,+1]}$ and define
\begin{equation}\Label{d3}
M_1^+=\underset p\bigcup I_p([1-\epsilon,1]).
\end{equation}
\bp\Label{p5s}
$M^+_1$ is a half space in a manifold $M_1$ of codimension $l-1$ with boundary $M$ and inward conormal $v'_o$ for
$v'_o$ close to $v_o$.
\ep
\bpf
We consider the mapping
\begin{equation}\Label{d4}
\Phi\colon \C^n\times\R^l\times[1-\epsilon,1]\to V',\,\,(w,x,r)\to 
I_p(r)\text{ for $p=(x+ih(x,w),w)$}.
\end{equation}
By Proposition~\ref{p4}, $\Phi$ is $C^{1,\beta'}$ in the complex 
of its arguments $(w,x)$  and $r$ up to $r=1$, and we have
$$
\Phi'_{(0,0,1)}(\C^n\times\R^l\times[1-\epsilon,1])=T_pM+\R^+v'_o.
$$
In particular $\Phi$ extends as a $C^{1,\beta'}$ mapping to 
$\C^n\times\R^l\times[1-\epsilon,1+\epsilon]$ whose image defines a 
manifold $M_1=\Phi(\C^n\times\R^l\times[1-\epsilon,1+\epsilon])$ which 
contains $M^+_1$ and satisfies  $T_pM^+_1=T_pM+\R v'_o$.

\epf

{\bf (e). End of proof of Theorems~\ref{t1.1}, \ref{t1.2}.}
First, we recall again that it suffices to prove Theorem~\ref{t1.2}. In fact, for $h_{I_j}=P_{I_j}+\mathcal O^{m_j+1}$ we have that $\langle\xi^o,P_{I_j}\rangle>0$ for $w_1\in\mathcal S$ implies $\langle\xi^o,h_{I_j}\rangle>0$ for $w_1\in\mathcal S$ and $|x_{I_i}|\leq c|w_1|^{m_i}$. Hence \eqref{1.5} is a consequence of \eqref{1.3}. Thus, let $f$ be a CR function on $M$. By the celebrated Baouendi-Treves approximation theorem of \cite{BT81}, $f$ is the uniform limit of polynomials on compact subsets of $M$. By the maximum principle it will extend to all analytic discs whose boundary is contained in this compact set. In particular it extends to the half-space $M^+_1$ of (b) for this is defined as the union of discs attached to $M$. On the other end it extends to a wedge $W$ with edge $M$ and directional cone, say $\Gamma$ by \cite{T90} since we are assuming that $M$ is of finite type. Thus by \cite{AH81} it will extend to a larger wedge $\hat W$ whose directional cone $\hat \Gamma$ is the convex hull of $\Gamma$ and $v'_o$ with $\langle\xi^o,v'_o\rangle>0$. In particular, for any $F$ holomorphic in $\hat W$, we have $\xi^o\notin WF(b(F))$. This completes the proof of Theorem~\ref{t1.2} and hence also Theorem~\ref{t1.1}.

We discuss some complements of our Theorems~\ref{t1.1}, \ref{t1.2}. 
We keep our choice of the $w_1$ - direction, select an index $i$, 
suppose $P_{I_{i}}=P_{I_{i}}(w_1)$ and \eqref{1.3}, or suppose \eqref{1.5}, and define
$\Gamma_{w_1,i}=\T{convex hull}\{v'_o\}$ where $v'_o$ ranges through the family of directions  produced by Th.~\ref{t1.1} or Th.~\ref{t1.2} for different directions $\xi_o$ and sectors $\mathcal S$. 
We use now the {\sc Ajrapetyan-Henkin} edge of the wedge theorem. 
In our setting it allows to state that all different directions of 
extension produced by Theorems~\ref{t1.1} or Theorem~\ref{t1.2}, and even those 
obtained as their convex combinations, can be collected to generate the
directional cone of a wedge of extension. Precisely,
for any $\epsilon$ there is a wedge $V'$ with edge $M$ and 
directional cone $\Gamma'_{w_1}$ satisfying 
$\Gamma'_{w_1}\subset(\Gamma_{w_1})_\epsilon$ and 
$\Gamma_{w_1}\subset(\Gamma'_{w_1})_\epsilon$ such that CR functions extend 
from $M$ to $V'$. (Here ${\cdot}_\epsilon$ denotes the $\epsilon$ 
conical neighborhood of $\cdot$. Also, in the above situation we will say that the cones $\Gamma_{w_1}$ and $\Gamma'_{w_1}$ are $\epsilon$-close.)
We can also play with different directions of the $w$-plane, say $w_k$. 
 Thus if we have equations of 
type $y_{w_k,I_i}=h_{w_k,I_i}$
with $h_{w_k,I_i}=\mathcal O^{m_i}$, then through Theorems~\ref{t1.1}, \ref{t1.2} we get directions $v'_{w_k,i}$ that we collect in a  cone 
\begin{equation}
\Gamma :=\sum_{k,i}\Gamma_{w_k,i}.
\end{equation} 
For this cone $\Gamma$ we have
\bp
\label{p3}
For any $\epsilon$ there is a wedge $V'$ with edge $M$ and 
directional cone $\Gamma'$ which is $\epsilon$-close to $\Gamma$, such that CR functions extend from 
$M$ to $V'$.
\ep
As already mentioned the proof is an immediate consequence of 
Theorem~\ref{t1.1} and \ref{t1.2} by the aid of the Ajrapetian-Henkin edge of the wedge theorem.
We want to discuss now about the dimension of $\Gamma_{w_k}$ and $\Gamma$.
Since we are dealing with various directions $w_k$'s, we will
write $m_{i,w_k}$, $l_{i,w_k}$ from now on. We have
\bp
\label{p1}
Assume that the equations  $y_{w_k,I_i}=h_{w_k,I_i}$ of $\tilde M=\C^l\times\C^1_{w_k}$ satisfy $h_{w_k,I_i}=\mathcal O^{m_i}$ and that $h_{w_k,I_i}$ is not $\tilde M$-harmonic. Then
$$
\T{dim}(\Gamma_{w_k})=\sum_{i}l_{w_k,i}.
$$
\ep
\bpf
We first prove that $\T{dim}(\Gamma_{w_k,1})=l_{w_k,1}$. We write $h_{w_k,1}=P_{w_k,1}(x,w_k)+\mathcal O^{m_{w_k},i}+1$ and know from the hypotheses that 
for any $\xi_o\in\R^{l_{w_k,1}}$, $\langle\xi_o,h_{w_k,I_1}(\tau 
w_k,x)\rangle$ is non-harmonic. In particular it is divisible 
by $|\tau|^2$ and hence, being of degree $m_{w_k,1}$ it has at most 
$2(m_{w_k,1}-2)$ zeroes on the unit circle $|\tau|=1$. In particular 
there is a sector of width $\geq\frac\pi{m_{w_k,1}-1}$ where it 
keeps constant sign and thus gives rise to a direction $v^o$ such 
that $\langle\xi_o,v^o\rangle\neq0$. If we play with all $\xi_o$ and 
all corresponding sectors, we conclude that these directions $v^o$ 
cannot be contained in any proper plane of $\R^{l_{w_k,1}}$.

We prove now the statement in full generality. For any $i$ we take a system of $l_{w_k,i}$ independent vectors 
$\xi\in\R^{l_{w_k,i}}$ and of corresponding sectors
$$
\mathcal 
S_\xi=\eta_ie^{i\theta_\xi}(1-\tau)^{\alpha_1}w_k,\quad\forall \tau\in\Delta,\,\,
\T{ with $\alpha_1$ satisfying $\frac1{m_{w_k,i}-1}>\alpha_1>\frac1{m_{w_k,i}}$}.
$$
We assume $\eta_1<<\eta_2\dots<<1$. This gives rise to a set of extension 
directions $v'=v'_{w_k,i,\xi,\mathcal S_\xi}$ of the type
$$
v'=(\eta_i^{m_1}v'_{I_1},\eta_i^{m_2}v'_{I_2},\dots,\eta_i^{m_i}v'_{I_i},\eta)\,\,\,\eta<<\eta_i\,\,\forall 
i,
$$
with the property that for each fixed $i$:
\begin{equation}
\T{dim}(\T{Span}_{\xi,\mathcal S_\xi}\{v'_{w_k,i,\xi,\mathcal 
S_\xi}\})=l_{w_1,i}.
\end{equation}
It is also clear, taking all $i$ and playing with 
different $\eta_i$, that
\begin{equation}
\T{dim}\left(\T{Span}_{i,\xi,\mathcal S_\xi}v'_{i,\xi,\mathcal 
S_\xi}\right)=
\sum_{i}l_{w_k,i}.
\end{equation}

\epf
Again, if we play with different directions $w_k$ we have the similar result as Proposition~\ref{p1} that is
\begin{equation}
\Label{collect}
\T{dim}\left(\sum_{k,i}\Gamma_{w_k,i}\right)=\sum_i\left(\T{rank}\{v'_{w_k,i}\}_{k}\right).
\end{equation}
(In this context the assumption that $M$ is of finite type that is $m_r<+\infty$ for a system of equations in Bloom-Graham normal form for the whole $M$, and not just for its $(l+1)$-dimensional sections $\tilde M$, precisely means on account of Proposition~\ref{p1}  and \eqref{collect} that $\T{dim}\Gamma=l$.)

\section{ H\"ormander's numbers of submanifolds of 
$\C^N$}
 Let $T^{1,0}M$ and $T^{0,1}M$ denote the bundles of vector fields 
tangent to 
$M$ which are holomorphic and antiholomorphic respectively. Let $T^\C 
M=TM\cap 
iTM$ be the complex tangent bundle to $M$; note that its 
complexification verifies 
$\C\otimes_\R T^\C M=T^{1,0}M\oplus T^{0,1}M$. Note that $\C\otimes_\R TM$ is integrable, 
that is closed under Lie brackets, but $\C\otimes_\R T^\C M$ is not, 
in general. We introduce a finite interpolation between 
$\C\otimes_\R T^\C M$ and $\C\otimes_\R TM$. We set
 $\mathcal L^1=\C\otimes_\R T^\C M$ and denote by $\mathcal{L}^j$ 
the 
distribution  of vector spaces spanned by
Lie brackets of holomorphic and antiholomorphic vector 
fields of length $\leq j$. Suppose that for an integer $m_1\geq 2$ we have
\begin{equation}\label{aI1}
\mathcal L^j_{p_o}=T^{1,0}_{p_o}M\oplus T^{0,1}_{p_o}M\,\,\,\forall 
j\leq 
m_1-1,\quad \mathcal L^{m_1}_{p_o}\underset\neq\supset  T^{0,1}_{p_o}M\oplus 
T^{1,0}_{p_o}M.
\end{equation}
Let $\text{dim}\frac{\mathcal L^{m_1}_{p_o}}{\mathcal L^1_{p_o}}=l_1$; in 
this situation it is usual to refer to $m_1$ as the 
{\em first H\"ormander number} of $M$ at $p_o$, and to $l_1$ as its 
{\em multiplicity}. 
In case $\mathcal L^j=\mathcal L^1$ for any $j$, we set $m_1=+\infty$ with 
multiplicity $l_1=l$. Next, we look for $m_2>m_1$ such that
\begin{equation}
\label{aI2}
\mathcal L^j_{p_o}=\mathcal L^{m_1}_{p_o}\,\,\forall j<m_2,\quad 
\mathcal L^{m_2}_{p_o}\neq\mathcal L^{m_1}_{p_o},
\end{equation}
and set $l_2=\T{dim}\left(\frac{\mathcal L^{m_2}_{p_o}}{\mathcal L^{m_1}_{p_o}}\right)$; 
again $m_2$ is possibly $+\infty$. We continue the above processus. We will call 
$M$ of finite type when commutators span the full $\C\otimes_\R T^\C_{p_o}M$. 
Thus the above chain will end with a  number $m_r<+\infty$ or $m_{r}=+\infty$ 
according to the case the type is finite or not.
We want to discuss now in greater detail about the first H\"ormander number. 
By the properties of 
linearity of commutators, one obtains easily the equivalence of 
\eqref{aI1}  to
\begin{multline}\Label{aI3}
[X_1,[X_2,\dots,[X_{j-1},X_j]\dots]\in T^{1,0}M\oplus T^{0,1}M 
\\
\forall X_i\in T^{1,0}M\oplus T^{0,1}M,\,\,\forall j\leq m_1-1
\end{multline}
\begin{multline}
\label{aI4}
[X_o^{\epsilon_1},[X_o^{\epsilon_2},\dots,[X_o,
\bar X_o]\dots]\notin 
T^{1,0}M\oplus T^{0,1}M
\\
\text{ for some $X_o$ and some choice of 
$X_o^{\epsilon_i}=X_o$ or $\bar X_o$}.
\end{multline}
One proves that commutators 
$[X_1[X_2,\dots,[X_{j-1},X_j]\dots]_{p_o}$, modulo 
$\C\otimes_\R T^\C M$
only depend on the initial values 
$X_1(p_o),\,X_2(p_o)\dots$ and not on the choice of the extended 
sections. This property is referred to as {\em tensoriality} of the 
iterated brackets of vector fields. 
We take a basis of equations $y_j=h_j,\,j=1,\dots,l$ for $M$ at $z_o=0$
with $h(0)=0$ and $\partial h(0)=0$ and also set $r_j=-y_j+h_j$ and $r=(r_j)$.  
We identify  $\frac{TM}{T^\C M}\overset\sim\to T_M\C^N$ by the complex structure $J$, and 
$T_M\C^N\overset\sim\to\R^l$ by the dual basis to $\partial r_j$. We 
look closely to $X_o$ in \eqref{aI4}, assume, say, 
$X_o(p_o)=w_o\partial_{w}$, and denote by $p'_o$ the projection of 
$p_o$ on the plane of $(x,w)$. 
We denote by $n-1$, resp. $m-1$ the occurrences of $X^{\epsilon_j}=X_o$ 
(resp. $X_o^{\epsilon_j}=\bar X_o$) in \eqref{aI4}. We can prove that 
\begin{equation}
\label{aI5}
\begin{cases}
\frac1{2i}[X_o^{\epsilon_1},\dots,X^{\epsilon_{j}},[X_o,\bar X_o],\dots]
(h)(p'_o)=0\quad\forall j< m_1-2
\\
\frac1{2i}[X_o^{\epsilon_1},\dots,X^{\epsilon_{m_1-2}},[X_o,\bar X_o],
\dots](h)(p'_o)=\partial_{w_o}^n\bar\partial_{w_o}^mh(p'_o).
\end{cases}
\end{equation}
This is a special case of subsequent Proposition~\ref{p6}. 
The above relation, together with the fact that harmonic terms can be 
removed by 
change of coordinates, makes \eqref{aI5} equivalent, in suitable 
coordinates, to
\begin{equation}\label{aI6}
\begin{cases}
\partial^\alpha_w\bar\partial_w^\beta h(p'_o)=0\,\,\forall 
|\alpha|+|\beta|\leq m_1-1 \\
\partial_w^\alpha h(p'_o)=0,\,\,\bar\partial_w^\alpha h(p'_o)=0,\,\,
\forall|\alpha|\leq m_1 \\
\partial_{w_o}^n\bar\partial_{w_o}^mh(p'_o)\neq 0 
\text{ for $X_o(p_o)=w_o\partial_{w}$ and for suitable $n+m=m_1$}.
\end{cases}
\end{equation}
We write also $\partial_{w_o}$ instead of $w_o\partial_w$ and consider  
the homogeneous term of lowest degree in the Taylor expansion of 
$h$ in the $w_o$-plane:
$$
g(\tau w_o)=\underset {\underset 
{m\geq1\,n\geq1}{m+n=k}}\sum\partial^m_{w_o}\bar\partial_{w_o}^nh(p'_o)\tau^m
\bar\tau^n.
$$
The above polynomial is real homogeneous and has some non-null 
coefficient on account of the third of \eqref{aI6}. Hence it has only a 
discrete set of zeroes for $|\tau|=1$
that is, for all $\theta\in [0,2\pi]$ but a discrete set, we have
$\sum\partial^m_{w_o}\bar\partial_{w_o}^nh(p'_o)e^{i(m-n)\theta}\neq0$. 
Sometimes we prefere to use the notation $\tilde w_o=e^{i\theta}w_o$ 
and then write in this notation
\begin{equation}\Label{notzero}
\underset {\underset {m\geq1\,n\geq1}{m+n=m_1}}\sum\partial^m_{\tilde 
w_o}\bar\partial_{\tilde w_o}^nh(p'_o)\neq0.
\end{equation}
We also denote by $v^o$ the vector in \eqref{notzero}. We remark that 
if $\xi_o\in\R^l$ verifies $\langle\xi_o,v^o\rangle\neq0$, then 
$$
\langle\xi_o,g( w_o)\rangle\gtrless0\text{ in a sector of the plane $\C_{w_o}$ of width 
$\geq\frac{\pi}{m_1-2}$.}
$$
In  fact each $g_i(\tau w_o)$ is divisible by $|\tau|^2$ and hence 
$|\tau|^{-2}\langle\xi,g(\tau w_o)\rangle$ has at most $m_1-2$ zeroes for 
$|\tau|=1$. Hence \eqref{aI5} or its equivalent version \eqref{aI6} imply 
our condition \eqref{1.2}.

To go further with our discussion, we need to fix  better our notations. 
We fix numbers $m_1<...<m_r$ (perhaps $m_r=+\infty$) and multiplicities 
$l_i$ with $\sum_il_i=l$. We take multiindices $I_1=(1,\dots,l_1),\dots,
I_r=(\sum_{i<r}l_i,\dots,l)$, give {\em weight} $m_i$ to the $x_{I_i}$ 
variables, and define the {\em weighted} vanishing order for a function 
$f=f(...x_{I_i}...,w)$ by putting $f=\mathcal O^{+\infty}$ when $m_r=+\infty$ 
and $f$ contains some monomial in the $x_{I_r}$'s, and, otherwise, putting 
$f=\mathcal O^m$ when $f(...t^{m_i}x_{I_i}...,tw)=\mathcal O(t^m)$. We then 
suppose that the equations of $M$ are presented according to increasing vanishing orders
\begin{equation}
\label{seminormal}
\begin{cases}
y_{I_1}=h_{I_1}
\\
...
\\
y_{I_r}=h_{I_r},
\end{cases}
\end{equation}
with $h_{I_i}=\mathcal O^{m_i}$ for any $i$. 
We point out that this is not necessarily the normal form in the {\sc Bloom Graham} sense. 
In fact we are not assuming that each $h_{I_i}$ is in the form $h_{I_i}=P_{I_i}(x_{I_1},
...,x_{I_{i-1}},w)$ with $\langle \xi,P_{I_i}\rangle$ non $M$-pluriharmonic for any $i$ 
and any $\xi\in\R^{l_i}$. (In this situation, {\em weighted homogeneity} does not serve any purpose.) 
To carry on our discussion, we need a description of a basis $\{X_j\}$ of vector fields 
for $T^{1,0}M$. We put $r_{I_i}=-y_{I_{i}}+h_{I_i}$, $r={}^t(r_1,\dots,r_{l})$, 
define an $(N-l)\times l$ matrix $A=(a_{jh})$ by
$$
A=-{}^t(\partial_w r)\,{}^t(\partial_z r)^{-1},
$$
and set $X_j=\sum_{h=1}^la_{jh}\partial_{z_h}+\partial_{w_j}$. We have
\begin{equation}
\label{i}
\sum_ha_{jh}\partial_{z_h}(r_{I_i})+\partial_{w_j}(r_{I_i})=0\,\,\forall i=1,\dots,r.
\end{equation}
Derivation of \eqref{i} yields
\begin{equation}
\label{ii}
\begin{cases}
\partial_{w\bar w}^\beta\partial_{x_{I_1}}^{\alpha_1}\dots\partial^{\alpha_{i-1}}_{x_{I_{i-1}}}(a_{j,I_i})=
0\T{ for $|\beta|+\sum_{j\leq i-1}m_j|\alpha_j|\leq m_i-2$}
\\
\sum_h \partial_{w\bar w}^\beta\partial_{x_{I_1}}^{\alpha_1}\dots
\partial^{\alpha_{i-1}}_{x_{I_{i-1}}}(a_{j,h})=-2i \partial_{w\bar w}^\beta
\partial_{x_{I_1}}^{\alpha_1}\dots\partial^{\alpha_{i-1}}_{x_{I_{i-1}}}\partial_{w_j}(r_{I_i})
\T{ for $|\beta|+\sum_{j\leq i-1}m_j|\alpha_j|\leq m_i-1$}.
\end{cases}
\end{equation}

Once the equations are ordered as in \eqref{seminormal}, we can introduce for any $i\leq r$ a diagram
\begin{equation}
\label{matrix}
\begin{matrix}
\frac{TM}{T^\C M}&\overset{\phi_1}\to &\R^l
\\
\downarrow& &\downarrow
\\
\frac{TM}{\mathcal L^{m_i-1}}&\overset{\phi_2}\to&\R^{l_i+...+l_r},
\end{matrix}
\end{equation}
where $\phi_1$ is defined by $[v]\mapsto J(v)(\partial r)$ and $\phi_2$ by 
$[v]\mapsto (Jv){}^t(\partial r_{I_i},...,\partial r_{I_r})$. 
We have to show that $\phi_2$ is well defined (in which case the diagram \eqref{matrix} is 
 commutative). 
To see this, we preliminarly remark that, just by the vanishing
condition in \eqref{seminormal}, we have $\{\partial r_{I_i},\dots,\partial r_{I_r}\}^{\perp_\C}=
\T{Span}_\R\{\partial_w,\bar\partial_w,\partial_{I_1},...,\partial_{I_{i-1}}\}$ (normal form being 
unessential for this conclusion). 
Thus our claim is a consequence of the following
\bp
\label{p5}
We have
$\mathcal L^{m_{i-1}}\subset\T{Span}\{\partial_w,\bar\partial_w,\partial_{x_{I_1}},...,\partial_{x_{I_{i-1}}}\}.$
\ep
\bpf
We have to show that
$$
[X_o^{\epsilon_1},...,[X_o^{\epsilon_{m_i-3}},[X_o,\bar X_o]...](r_{I_{j}})
(p_o)=0\,\,\forall j\geq i\T{ and for any $\epsilon$}.
$$
We recall \eqref{i} and \eqref{ii} and fix $j=i$. We use the notation $[\cdot,\cdot]^k$ 
to denote brackets of $X_o$ or $\bar X_o$ performed $k-1$ times. We assume, for instance, 
$X_o^{\epsilon_1}=\sum_h a_{1h}\partial_{z_h}+\partial_{w_1}$, and begin by remarking that 
\begin{equation}
\label{iii}
\begin{split}
[\cdot,\cdot]^{m_i-1}(r_{I_i})&=[\sum_h a_{1h}\partial_{z_h}+\partial_{w_1},[\cdot,\cdot]^{m_i-2}](r_{I_i})
\\
&=[\partial_{w_1},[\cdot,\cdot]^{m_i-2}](r_{I_i}),
\end{split}
\end{equation}
due to $a_{1h}(p_o)=0$ and $[\cdot,\cdot]^{m_i-2}(a_{1h})=0$.
Continuing in this way we end up with
\begin{multline}
\label{iv}
[\partial_{w_1\bar w_1}^{\beta},[\sum_ha_{1h}\partial_{z_h}+\partial_{w_1},
\sum_h\bar a_{1h}\bar\partial_{z_h}+\bar\partial_{w_1}]](r_{I_i})
\\
=\partial_{w_1\bar w_1}^\beta\left(\sum_ha_{1h}\partial_{z_h}(\bar a_{1I_i})-\frac i2\partial_{w_1}
(\bar a_{1I_i})+\sum_ha_{1h}\partial_{z_h}\bar\partial_{w_1}(r_{I_i})+\partial_{w_1}
\bar\partial_{w_1}(r_{I_i})\right)+\dots
\end{multline}
where  
$\beta$ is a biindex of length $|\beta|=m_i-3$ and
the dots denote similar terms as  the four in the right hand side of \eqref{iv}. Now:
\begin{gather}
\label{v}
\partial_{w_1\bar w_1}^\beta\partial_{w_1}(\bar a_{1I_i})=0\T{ (by \eqref{iii})},
\\
\label{vi} 
\partial_{w_1,\bar w_1}^\beta\left(\sum_ha_{1h}\partial_{z_h}\bar\partial_{w_1}(r_{I_i})\right)=
\sum_{\gamma+\delta=\beta}\sum_h\partial_{w_1,\bar w_1}^\gamma(a_{1h})\partial^\delta_{w_1,
\bar w_1}\partial_{z_h}\bar\partial_{w_1}(r_{I_i}).
\end{gather}
Thus, if $h\in I_{j}$ for $j\geq i$, the above term is clearly $0$. Otherwise, 
either $|\gamma|\leq m_h-2$ and hence $\partial_{w_1\bar w_1}^\gamma(a_{1h})=0$,  
or else $|\delta|\leq m_i-2-m_h$ and hence $\partial_{w_1\bar w_1}^\delta
\partial_{z_h}\bar\partial_{w_1}(r_{I_i})=0$. By the same reason, we have 
for the remaining term in \eqref{iv}: $\partial_{w_1\bar w_1}^\beta
\left(\sum_ha_{1h}\partial_{z_h}\right)(\bar a_{1I_i})=0$. Finally, 
$\partial_{w_1\bar w_1}^\beta\partial_{w_1}\bar\partial_{w_1}(r_{I_i})$ is also $0$ again by \eqref{iii}.
The proof is complete.

\epf
\br
Note that $\phi_2$ is an isomorphism precisely when we have in fact equality in Proposition~\ref{p5}. 
But this is  equivalent as to asking that the equations~\eqref{seminormal} are in normal form. 
\er
Let us choose a vector field $X_o\in T^{1,0}M$ with $X_o(p_o)=w_o\partial_w$; 
we will also use the notation $\partial_{w_o}$ instead of $w_o\partial_w$.  We have
\bp
\label{p6}
\begin{equation}
\label{aI7}
[X_o^{\epsilon_1},...,[X_o^{\epsilon_{m_i-2}},[X_o,\bar X_o]...](r_{I_j})
(p_o)=0\,\,\forall j>i \T{ and any $\epsilon$}. 
\end{equation}
If moreover $\langle\xi_o,h_{I_i}\rangle$, restricted to $\C_{w_o}\times\R^l_x$, is in the form
$P+\mathcal O^{m_i+1}$ for $P=P(w_o)$ homogeneous of degree $m_i$ with $m_i<+\infty$,  then
\begin{equation}
J[X_o^{\epsilon_1},...[X_o^{\epsilon_{m_1-2}},[X_o,\bar X_o]...]\langle\xi_o,r_{I_i}
\rangle(p_o)=-2\partial^n_{w_o}\bar\partial_{w_o}^m(\langle\xi_o,h_{I_i}\rangle)(p'_o).
\end{equation}
\ep
\bpf
The first statement is a variant of Proposition~\ref{p5}. As for the second,
 in the same way as in the proof of Proposition~\ref{p5}, we get for a suitable $|\beta|=m_i-2$
\begin{equation}
\label{vii}
\begin{split}
[X_o^{\epsilon_1},...,&[X_o^{\epsilon_{m_1-2}},[X_o,\bar X_o]...](r_{I_i})=
[\partial^\beta_{w_1\bar w_1},[\sum_ha_{1h}\partial_{z_h}+\partial_{w_1},
\sum_h \bar a_{1h}\bar\partial_{z_h}+\bar\partial_{w_1}]](r_{I_i})
\\
&=\partial_{w_1\bar w_1}^\beta\left(\sum_ha_{1h}\partial_{z_h}(\bar a_{I_i})-
\frac i2\partial_{w_1}(\bar a_{I_i})+\sum_ha_{1h}\partial_{z_h}\bar\partial_{w_1}(r_{I_i})+
\partial_{w_1}\bar\partial_{w_1}(r_{I_i})\right)+\dots,
\end{split}
\end{equation}
where the dots denote similar terms. Now the fourth term disappears by elimination with the 
terms in the dots (where it appears with opposite sign).
 The first and third term are not 0, in general. However, they vanish if we apply vector 
 fields not to the whole $r_{I_i}$ but just to $\langle\xi_o,r_{I_i}\rangle$ on account 
 of the hypothesis of {\em semirigidity} contained in the second statement of the proposition. 
 Thus, for the third term, we have
$$
\partial_{w_1\bar 
w_1}^\beta\left(\sum_ha_{1h}\partial_{z_h}\bar\partial_{w_1}\right)\langle\xi_o,
r_{I_i}\rangle=\sum_{\gamma+\delta=\beta}\partial_{w_1\bar w_1}^\gamma(\sum_ha_{1h}
\partial^\delta_{w_1\bar w_1}\partial_{z_h}\bar\partial_{z_h}\bar\partial_{w_1}\langle\xi_o,r_{I_i}\rangle.
$$
 Again, if $|\gamma|\leq m_h-2$ then $\partial_{w_1\bar w_1}^\gamma a_{ih}=0$. 
 If, instead, $|\delta|\leq m_i-1-m_h$ then $\partial^\delta_{w_1\bar w_1}
 \partial_{z_h}\bar\partial_{w_1}\langle\xi_o,r_{I_i}\rangle=
 \partial^\delta_{w_1\bar w_1}\partial_{z_h}\bar\partial_{w_1}
 (P+\mathcal O^{m_i+1})=0$. In the same way one proves that the 
 first term in the second line of \eqref{vii} is 0.
The only term which survives is the seond (which also appears, with 
the same sign in the dots terms). We have thus got
\begin{equation}
\begin{split}
J[X_o^{\epsilon_1},...[X_o^{\epsilon_{m_1-2}},[X_o,\bar X_o]...]
\langle\xi_o,r_{I_i}\rangle&=-\frac i2(\partial^\beta_{w_1\bar w_1}
\partial_{w_1}\langle\xi_o,\bar a_{1I_i}\rangle+\dots
\\
&=-\frac i2(\partial^\beta_{w_1\bar w_1}\partial_{w_1}\bar\partial_{w_1}
\langle\xi_o,r_{I_i}\rangle\frac2i+\dots)
\\
&=-2 \partial^\beta_{w_1\bar w_1}\partial_{w_1}\bar\partial_{w_1}\langle\xi_o,r_{I_i}\rangle.
\end{split}
\end{equation}
This completes the proof of the proposition.

\epf
We assume now that for some vector field $X_o$ with $X_o(p_o)=\partial_{w_o}$, for some 
$\epsilon=(\epsilon_1,...,\epsilon_{m_i-2})$, and for some $\xi_o\in\R^{l_i}$, we have
\begin{equation}
\label{aI9}
[X_o^{\epsilon_1},...,[X^{\epsilon_{m_i-2}},[X_o,\bar X_o]...]\notin\T{Span}\{\partial_w,
\partial_{\bar w},\partial_{x_{I_1}},\cdot,\cdot,\partial_{x_{I_{i-1}}}\}
\end{equation}
and
\begin{equation}
\label{aI10}
\langle\xi_o,h_{I_i}\rangle|_{\C_{w_o}\times\R^l_x}=P(w_o)+\mathcal O^{m_i+1}.
\end{equation}
It follows that $P(w_o)=|w_o|^2Q(w_o)$ with $Q$ real homogeneous of degree $m_i-2$. 
Since  $Q$ has at most $m_i-2$ zeroes on the circle $|w_o|=1$, then
\begin{equation}
\label{aI11}
 P\gtrless0\T{ for $w_o$ in a sector of width }>\frac\pi{m_i-2}.
\end{equation}
Hence we enter in the hypotheses of Theorem~\ref{t1.2} and conclude that CR 
functions on $M$ extend to a new direction $v^o$ satisfying $\langle\xi_o,v^o\rangle\gtrless0$. 
Note that in that Theorem normal equations as in \eqref{1.3} are not needed. 
What is really needed is, for equations as \eqref{seminormal}, to assume 
$\langle\xi_o,h_{I_i}\rangle=P(w_o)+\mathcal O^{m_i+1}$ and $P\geq0$ 
(or $P\leq0$) in a sector $>\frac\pi{m_i}$.

Naturally, if the equations are normal, we have the significant simplification that 
$\mathcal L^{m_i}=\T{Span}\{\partial_{w},\partial_{\bar w},\partial_{x_{I_i}},\dots,
\partial_{I_{m_{i-1}}}\}$. Thus vector fields $X_o$ which satisfy   \eqref{aI9} do exist. 
If for one of them, with, say, $X_o(p_o)=\partial_{w_o}$, and for some $\xi_o\in\R^{l_i}$,  
\eqref{aI10} is also satisfied, then Proposition~\ref{p3} yields CR extension to some $v^o$ 
with $\langle\xi_o,v^o\rangle\gtrless0$. 

\section{  Comparison with {\sc Boggess-Pitts} 
\cite{BP85}}
Let $M$ be a manifold of class $C^{k+2}$  which satisfies \eqref{1.1}  
with $g$ homogeneous of degree $k$ and non $M$-harmonic
(in particular whose first H\"ormander number is $m_1=k$. Remember that in this situation \eqref{1.2} 
is also satisfied.
Let $v$ be the direction normal to $M$ given by the formula
\begin{equation}
\label{cepsilon}
 v=\underset\epsilon\sum 
C_\epsilon[X_o^{\epsilon_1},X_o^{\epsilon_2},\dots,[X_o,\bar 
X_o]\dots](r)(p_o)\text{ where $C_\epsilon:=\frac{1}{\epsilon^+!\epsilon^-!}$}
\end{equation}
with $\epsilon^+$ and $\epsilon^-$ denoting the occurences $X_o^{\epsilon_i}=X_o$ and 
$X_o^{\epsilon_i}=\bar X_o$ respectively. Note that the last two occurences are fixed as 
$X_o^{\epsilon_{k-1}}=X_o$ and $X_o^{\epsilon_k}=\bar X_o$. Let 
$X_o(p_o)=\partial_{w_o}$. By tensoriality of brackets and by the 
combinatorial remark that the number of choices of $\epsilon$'s which 
give rise to the same pair of occurences $m$, $n$ is 
${{k-2}\choose{m-1}}$ one gets
\begin{equation}\Label{v}
v=\underset{\underset{m\geq1\,n\geq1}{m+n=k}}
\sum{{k-2}\choose{m-1}}\frac1{m!n!}\partial_{w_o}^m\bar\partial_{w_o}^nh(p'_o).
\end{equation}
Again,  once the complex plane of $X_o(p_o)$ is fixed, in our case 
the $w_o$-plane, there might be many vectors $v=v_\phi$ produced 
through \eqref{v} just by replacing $w_o$ by $e^{i\phi} w_o$. The result 
by {\sc Boggess-Pitts} \cite{BP85} is that for each of these vectors 
$v$, one obtains CR extension from $M$ to $M'$ where $M'$ points to a 
direction $v'$ close to $v$. We first discuss this extension 
in case $M$ is a hypersurface of $\C^N$ defined, in coordinates $(z,w)\in\C^1\times\C^n,\,w=(w_1,w')$, for a 
pair of even integers $k$ and $p$  with $p\leq k-2$, for a choice of a coefficient $a\geq0$, 
and with the notation  $w_o=(1,0,\dots)$,
by an equation
\begin{equation}
\Label{vi}
y_1=|w_1|^k+a|w|^{k-p}\Re w_1^p+
\left(O(|x_1|^2+|w_1|^{k+1}+|x_1||w_1|+|w||w'|\right).
\end{equation}
We denote by $g=g(w_1)$ the homogeneous polynomial in the right side of \eqref{vi}. 
With $p_o=0$ and $X_o=\partial_{w_1}$  and with the notation $k-2=p+2q$, we have extension in directions 
$v_\phi=(ic_\phi,0,\dots)$
for 
$$
c_\phi={{k-2}\choose{\frac k2-1}}+a\cos(p\phi){{k-2}\choose{p+q}}.
$$
In particular if we look for  extension {\em down}, that is for $v_\phi$ with negative first component, 
we have to require $k\geq 4$, $p\geq2$. Then  $v_\phi<0$ will occur exactly for $\phi=\frac{\pi}{p}$ 
(which yields $\cos(p\phi)=-1$) and  
$$
a\geq\frac{(p+q)!q!}{(\frac k2-1)!(k-1-\frac k2)!}.
$$ 
We compare the above condition with that which is given by sector property.
We consider the restriction of $g$ on the unit circle $w_1=e^{i\theta}$ 
given by $g(e^{i\theta})=1+a\cos(p\theta)$. It is clear that for any choice of $a$ we have $g\geq0$ in a sector 
of width bigger than $\frac\pi p$ which is in turn bigger than $\frac\pi k$. Hence by Theorem 1.1 we get 
holomorphic extension {\em up}.

If we search, instead, for extension {\em down}, we can use the following result which generalizes similar 
conclusions by {\sc Baouendi-Treves} \cite{BT82} concerning the case $k=4$.
\bp
\Label{p7}
We have
\begin{equation}
g<0\T{ in a sector of width $>\frac\pi{k}$,}
\end{equation}
if and only if
\begin{equation}
a>\frac1{\cos\left(\frac{p\pi}{2k}\right)}.
\end{equation}
\ep
\bpf
Let $a>0$;
it is clear that $1+a\cos(p\theta)$ attains its minimum at $\theta=\frac\pi{p}$. It is also clear that in 
order that the sector where $g<0$ has angle bigger than $\frac\pi{k}$ it is necessary and sufficient that 
$$
a\cos(\pi+\frac{p\pi}{2k})<-1,
$$
which is equivalent to the condition in the statement of the proposition.

\epf
We also have the following statement which shows 
necessity of sector property for
holomorphic extendibility.
\bp
\Label{p8}
Let $p$ divide $k$ and  $a\leq \frac1{\cos\left(\frac{p\pi}{2k}\right)}$. Then for $b=\frac pk\T{tg}\left(\frac{p\pi}{2k}\right)$, 
we have that the trigonometric polynomial
$g_1=1-a\cos(p\theta)+b\cos(k\theta)$  verifies
$g_1\geq0\,\,\forall \theta$. 
In particular,
if in addition the plane of the $w$ variables has dimension 1,
by adding another harmonic term $\epsilon \T{sin}({k\theta})$ we can achieve $g_1(w_1)\geq c_1|w_1|^m$ for $c_1>0$.
\ep
\bpf
Since $a\leq \frac1{\cos\left(\frac{p\pi}{2k}\right)}$, then for $g_1\geq0$ it will suffice:
\begin{equation}
\Label{sector}
b\cos(k\theta)\geq\frac1{\cos\left(\frac{p\pi}{2k}\right)}\cos(p\theta)-1.
\end{equation}
It is clear that it will suffice to take $b$ such that
\begin{equation}
\Label{bi}
(b\cos(k\theta)'|_{-\frac{\pi}{2k}}=\frac 1 {\cos\left(\frac{p\pi}{2k}\right)}\cos
(p\theta)'|_{-\frac{\pi}{2k}}.
\end{equation}
In fact this choice of $b$ will imply that the derivative on the left of \eqref{bi} dominates (respectively is 
dominated by) the one on the right in the interval
$[-\frac\pi{2k},0]$ (respectively in $[-\frac\pi k,-\frac\pi{2k}]$). Hence \eqref{sector} holds in the  interval 
$[-\frac\pi{k},0]$
and also, by symmetry, in the whole $[-\frac\pi{k},+\frac\pi{k}]$. 
It also holds trivially in the remaining part of $[-\frac\pi{p},+\frac\pi{p}]$.
On the other hand this is a complete cycle of the trigonometric 
function $1-a\T{cos}(p\theta)+b\T{cos}(k\theta)$ due to the 
assumption that $p$ divides $k$.
\epf
\bc
\Label{c4.1}
Let $M$ be a hypersurface in $\C^N$ defined by \eqref{vi} and assume that $p$ divides $k$. If $a\leq \frac1{\cos\left(\frac{p\pi}{2k}\right)}$, then there are CR functions $f\in CR_M$ which do not extend {\it down}. 
\ec
\bpf
In new complex coordinates we can arrange that $M\subset\{y_1>0\}$. Since $\{y_1>0\}$ is pseudoconvex, the conclusion follows.
\epf
The comparison between the conditions related to \cite{BP85} and to sector property is expressed by 
\bl
\Label{l10}
Let $k-2=p+2q$. Then
\begin{equation}
\Label{compare}
\frac{(p+q)!q!}{(\frac k2-1)!(\frac k2-1)!}>\frac1{\cos\left(\frac {p\pi}{2k}\right)}.
\end{equation}
\el
\bpf
The most delicate case is when $p=2$. In this case \eqref{compare} becomes
$$
\frac{\frac k2!(\frac k2-2)!}{(k-1)!(\frac k2-1)!}>\frac 1{\cos\left(\frac\pi k\right)},
$$
or else
$$
\frac{\frac k2}{\frac k2-1}>\frac1{\cos\left(\frac\pi k\right)}.
$$

\epf
Hence the method of sectors  is sharper than that of \cite{BP85}.
 In particular it yields extension {\em down} for an extra range of values of $a$ that is for 
 $\frac1{\cos\left(\frac {p\pi}{2k}\right)}\leq a<\frac{(p+q)!q!}{(\frac k2-1)!(\frac k2-1)!}$. 
 The above conclusions are generalizations of former results by {\sc Baouendi-Treves} \cite{BT82}. 

We pass now to the higher codimensional case. We discuss CR-extension 
for 
$M\subset\C^3$ defined in coordinates $(z_1,z_2,w)$ by the system
\begin{equation}
\Label{example}
\begin{cases}
y_1=|w|^k+a|w|^2\Re w^p+O(|x|^2+|w|^{k+1}+|x||w|)
\\
y_2=|w|^k+O(|x|^2+|w|^{k+1}+|x||w|)
\end{cases}
\end{equation}
We also denote by $g=(g_j)_j, \,\,j=1,2$ the vector with polynomial 
entries on the right of \eqref{example} and, for $\xi\in\R^2$, we use the notation 
$g_\xi=\langle \xi,g\rangle$. We
can express the extension directions $v_\phi$ by \cite{BP85} as
\begin{equation}
v_\phi=({{k-2}\choose{\frac k2-1}}+a\cos(\phi){{k-2}\choose{p+q}},  {{k-2}\choose{\frac k2-1}})^t,
\end{equation}
where $\cdot^t$ denotes transposition. Let us search for $v_\phi$ whose first component is $<0$. 
The first occurence, which takes place for $\phi=\pi$ is when $a>
\frac{(p+q)!q!}{(\frac k2-1)!(\frac k2-1)!}$. 
 In this case extension 
to directions arbitrarily close to 
$v=({{k-2}\choose{\frac k2-1}}-a{{k-2}\choose{p+q}},  {{k-2}\choose{\frac k2-1}})^t$  holds according to 
\cite{BP85}. If we look, instead, to our sector property and search for $v$ 
whose first component is $<0$ and the second is, say, $>0$,
we are lead to the sector property of $g_\xi$ for suitable $\xi$ with 
$\xi_1<0$ and $\xi_2<0$. The condition reads in this case
%
\begin{equation}
\begin{split}
g_\xi(\theta)&= \xi_1(1+a\cos(p\theta))+\xi_2
\\
&>0\text{ in a sector of width $>\frac\pi k$},
\end{split}
\end{equation}
that is
$$
1+a\frac{\xi_1}{\xi_1+\xi_2}\cos(p\theta)<0
\text{ in a sector 
with angle $>\frac\pi k$}.
$$
We write $a_\xi=a\frac{\xi_1}{\xi_1+\xi_2}$. Then the sector where 
$g_\xi>0$ is centered at  $\theta=\frac\pi p$  and its width 
is $>\frac\pi k$ if and only if 
$$
a_\xi>\frac 1{\cos\left(\frac{p\pi}{2k}\right)}.
$$
Now the first such an occurence is for $a>\frac1{\cos\left(\frac{p\pi}{2k}\right)}$ and 
for $\xi$ close to $(-1,1-a{\cos\left(\frac{p\pi}{2k}\right)})^t$. Hence we get 
extension to vectors 
$v$ with $\langle v,\xi\rangle>0$ according to Theorem 1.1. Direct inspection of the second equation of $M$ 
shows that $v_2>0$. Also, extension to directions of a conic neighborhood of the diagonal is evident. In conclusion, 
using also the edge of the wedge theorem, we get extension to all intermediate directions, among whose some which is 
close to 
 $v=(-a{\cos\left(\frac{p\pi}{2k}\right)} +1,1)^t$. 
We write now $a_1=a\frac{(\frac k2-1)!(k-1-\frac k2)!}{(p+q)!q!}$ and $a_2=a{\cos\left(\frac{p\pi}{2k}\right)}$. 
According to Lemma~\ref{l10} we always have $a_1<a_2$. (In the simplest cases we have $\frac{a_2}{a_1}=
\frac2{\sqrt2}$ when $k=4,\,\,p=2$ and $\frac {a_2}{a_1}=\frac{3}{\sqrt3}$ for $k=6,\,\,p=4$.)
Summarizing up we get
{\bf (1)} Extension for an extra range of values of $a$ that is 
$a_1<a\leq 
a_2$ which were not taken care of by \cite{BP85}. For this purpose the higher codimension is not really needed; 
the examples by Rea and Baouendi-Treves would suffice as well.
{\bf (2)} Extension to a  wedge $V'$ with bigger directional cone 
$\Gamma'$ even 
for common values of $a>a_1$. (Here codimension $>1$ is really essential.)
In fact in \cite{BP85} the cone is 
$$
\Gamma=\{(y_1,y_2):y_1>-|y_2|(a_1-1)\},
$$
whereas in our case it is
$$
\Gamma=\{(y_1,y_2):\,y_1>-|y_2|( a_2-1).
$$


\begin{thebibliography}{BER99}
\bibitem[1]{AH81} {\bf Ajrapetyan,~R.A.; Henkin,~G.M.} --- Analytic 
continuation of CR-functions through the "edge of the 
wedge".{\em  Sov. Math., Dokl}, (1981),  129--132
\bibitem[2]{B91} {\bf  Boggess,~A.} --- CR manifolds and the 
tangential 
Cauchy-Riemann complex. {\em Studies in Adv. Math. CRC Press}, (1991).
\bibitem[3]{BER99}{\bf  Baouendi,~M.S.; Ebenfelt,~P.; 
Rothschild,~L.P.} --- Real Submanifolds in Complex Space and Their 
Mappings.
{\em Princeton Math. Series, Princeton Univ. Press}, (1999).
\bibitem[4]{BT81} {\bf Baouendi,~M.S.; Treves,~F.} ---
 A property of the functions and distributions annihilated by a 
locally integrable system of complex vector fields.
{\em Ann.~Math.}
{\bf  114} (1981),
 387--421.
\bibitem[5]{BR87}{\bf Baouendi,~M.S; Rothschild,~L.} --- 
Normal forms for generic manifolds and holomorphic extension of CR functions. 
{\em J. Diff. Geom.}, {\bf 25} (1987), 431--467.
\bibitem[6]{BG77} {\bf Bloom,~T., Graham,~I.} --- On type conditions 
for generic real submanifolds of $\C^n$. {\em Invent. Math.}, {\bf 
40} (1977), 217--243.
\bibitem[7]{BZZ02}{\bf  Baracco,~L.; Zaitsev,~D., Zampieri,~G} --- 
CR extension from wedges on manifolds of higher type.
(2002), to appear
\bibitem[8]{BP85} {\bf Boggess,~A.; Pitts,~J.}
--- CR extension near a point of higher type.
{\em Duke Math. J.},
{\bf  52 (1)} (1985), 67--102.
\bibitem[9] {BP82}
{\bf Boggess,~A.; Polking,~J.C.} ---
Holomorphic extension of CR functions.
{\em Duke Math.~J.}
{\bf  49} (1982), 757--784.
\bibitem[10]{EG99}
{\bf  Eastwood,~M.C.; Graham,~C.R.}
--- An Edge-of-the Wedge Theorem for Hypersurface CR Functions
{\em  J. Geom. Anal.} {\bf 11 (4)} (2001), 589--602.
\bibitem[11]{Tr85}
{\bf Trepreau,~J.M.}
--- Sur le prolongement holomorphe des fonctions CR definies sur une hypersurface r\'eelle de classe $C^2$ dans $\C^n$. 
{\em Invent. Math.},
{\bf 83 }
(1986),
583--592.
\bibitem[12]{T90}
{\bf Tumanov,~A.E.}
--- Extension of CR-functions into a wedge.
{\em Mat. Sb.},
{\bf 181 (7)}
(1990),
951--964.
\bibitem[13]{T95}
{\bf Tumanov, A.E.}
--- Extending CR functions from manifolds with boundaries.
{\em Math. Res. Lett.}
{\bf 2 (5)}
(1995),
629--642. 
\bibitem[14]{T96}
{\bf Tumanov,~A.E.}
--- Analytic discs and the extendibility of CR functions. {\em 
Integral geometry, Radon transforms and complex analysis 
(Venice, 1996) Lecture Notes in Math., Springer, Berlin}, {\bf 1684} 
(1998),  123--141.
\bibitem[15]{ZZ01}
{\bf Zaitsev,~D; Zampieri,~G.}
--- Extension of CR-functions into weighted wedges through 
families of nonsmooth analytic discs.
{\em Transactions of the AMS} (2003).
\bibitem[16]{ZZ02}
{\bf Zaitsev,~D; Zampieri,~G.} --- Extension of CR functions on 
wedges.
{\em Math.~Ann.}, (2003).
\bibitem[17]{BT82}
{\bf Baouendi,~M.S.; Treves,~F.} --- About holomorphic extension of CR 
functions on real hypersurfaces in complex space.
{\em Duke Math. J.}
{\bf 51}
(1984),
77--107.


\end{thebibliography}
\end{document}